\documentclass[10pt]{article}

\usepackage{amscd,amsmath, amssymb, fancyhdr}

\usepackage[backref=page]{hyperref}
\renewcommand*{\backrefalt}[4]{%
	\ifcase #1 (Not cited.)%
	\or        (Cited on page~#2.)%
	\else      (Cited on pages~#2.)%
	\fi}

\hypersetup{
	colorlinks   = true,
	citecolor    = magenta,
	linkcolor    = blue,
	urlcolor     = magenta	
}

\numberwithin{equation}{section}

\newcommand{\version}{version 1.2,\ \ September 15, 2025}

\def\eqref#1{(\ref{#1})}
\newcommand{\goth}{\mathfrak}

\newcommand{\arrow}{{\:\longrightarrow\:}}

\def\C{{\Bbb C}}

\def\1{\sqrt{-1}\:}
\newcommand{\restrict}[1]{{\left|_{{\phantom{|}\!\!}_{#1}}\right.}}
\newcommand{\cntrct}                % contraction with a vector field
{\hspace{2pt}\raisebox{1pt}{\text{$\lrcorner$}}\hspace{2pt}}

\newcommand{\calo}{{\cal O}}

\renewcommand{\bar}{\overline}
\renewcommand{\phi}{\varphi}
\renewcommand{\epsilon}{\varepsilon}
\renewcommand{\geq}{\geqslant}
\renewcommand{\leq}{\leqslant}

\newcommand{\Discr}{\operatorname{{\sf Discr}}}

\newcommand{\codim}{\operatorname{codim}}

\newcommand{\rk}{\operatorname{rk}}

\newcommand{\Spec}{\operatorname{Spec}}

\newcounter{Mycounter}[section]
\newcounter{lemma}[section]
\newcounter{claim}[section]
\renewcommand{\theclaim}{{Claim \thesection.\arabic{claim}}}
\newcommand{\claim}{%
    \setcounter{claim}{\value{Mycounter}}
    \refstepcounter{claim}
    \stepcounter{Mycounter}
    {\noindent \bf \theclaim:\ }}

\newcounter{sublemma}[section]
\newcounter{corollary}[section]
\renewcommand{\thecorollary}{{Corollary \thesection.\arabic{corollary}}}
\newcommand{\corollary}{%
    \setcounter{corollary}{\value{Mycounter}}
    \refstepcounter{corollary}
    \stepcounter{Mycounter}
    {\noindent \bf \thecorollary:\ }}

\newcounter{theorem}[section]
\renewcommand{\thetheorem}{{Theorem \thesection.\arabic{theorem}}}
\newcommand{\theorem}{%
    \setcounter{theorem}{\value{Mycounter}}
    \refstepcounter{theorem}
    \stepcounter{Mycounter}
    {\noindent \bf \thetheorem:\ }}

\newcounter{conjecture}[section]
\newcounter{proposition}[section]
\renewcommand{\theproposition}
      {{Proposition \thesection.\arabic{proposition}}}
\newcommand{\proposition}{%
    \setcounter{proposition}{\value{Mycounter}}
    \refstepcounter{proposition}
    \stepcounter{Mycounter}
    {\noindent \bf \theproposition:\ }}

\newcounter{definition}[section]
\renewcommand{\thedefinition}
      {{Definition~\thesection.\arabic{definition}}}
\newcommand{\definition}{%
    \setcounter{definition}{\value{Mycounter}}
    \refstepcounter{definition}
    \stepcounter{Mycounter}
    {\noindent \bf \thedefinition:\ }}

\newcounter{example}[section]
\newcounter{remark}[section]
\renewcommand{\theremark}{{Remark \thesection.\arabic{remark}}}
\newcommand{\remark}{%
    \setcounter{remark}{\value{Mycounter}}
    \refstepcounter{remark}
    \stepcounter{Mycounter}
    {\noindent \bf \theremark:\ }}

\newcounter{problem}[section]
\newcounter{question}[section]
\newcommand{\pstep}{{\bf Proof. Step 1:\ }}
\newcommand{\proof}{{\bf Proof:\ }}

\makeatletter

\@addtoreset{equation}{section} \@addtoreset{footnote}{section}
\makeatother

\def\blacksquare{\hbox{\vrule width 5pt height 5pt depth 0pt}}
\def\endproof{\blacksquare}

\begin{document}

%%%%%%%%%%%%%%%%%%%%%%%%%%%%%%%%%%%%%%%%%%%%%%%%%%%%%%%%%%%%
\begin{center}
{\LARGE\bf
Torsor N\'eron models of hyperk\"ahler manifolds \\[4mm]
}
%%%%%%%%%%%%%%%%%%%%%%%%%%%%%%%%%%%%%%%%%%%%%%%%%%%%%%%%%%%%%

Ljudmila Kamenova\footnote{Partially supported 
by award SFI-MPS-TSM-00013537 from the Simons Foundation International},  
Misha Verbitsky\footnote{Partially supported 
by FAPERJ SEI-260003/000410/2023 and CNPq - Process 310952/2
021-2. 

 }

\end{center}

{\small \hspace{0.10\linewidth}
\begin{minipage}[t]{0.85\linewidth}
{\bf Abstract.}  Let $X$ be a compact hyperk\"ahler 
manifold equipped with a Lagrangian fibration
$\pi:\; X \to Y$, and $X'$ the smooth locus of $\pi$.
We prove that over a complement to a codimension $\geq 2$ subset in $Y$,
the projection $\pi:\; X' \to Y$ has a natural structure of
a torsor over an abelian group bundle, which can be understood
as a complex analytic variant of the N\'eron model construction.
This gives an independent proof of a result by Y.-J. Kim.
\end{minipage}
}

\tableofcontents

%%%%%%%%%%%%%%%%%%%%%%%%%%%%%%%%%%%%%%%%%%%%%%%%%%%%%%%%%%%%
\section{Introduction}
%%%%%%%%%%%%%%%%%%%%%%%%%%%%%%%%%%%%%%%%%%%%%%%%%%%%%%%%%%%%

The results of this paper were originally suggested
by F. Bogomolov (see also \cite{_BHPT:abelian_CY_}),
and an early version of our proof appeared in
our paper \cite{_BKV:Sections_}, joint with Bogomolov.
After \cite{_BKV:Sections_} appeared on the arXiv, 
we had fruitful discussions with J. Koll\'ar, who
posted the paper \cite{_Kollar:Neron_} which
helped to significantly simplify some of our arguments.

After a while, we realized that our results
can be applied in a more general situation, and
might have applications outside of the scope of
\cite{_BKV:Sections_}. 

The original idea of Bogomolov was to apply the
construction of N\'eron models to produce smooth
sections of Lagrangian fibrations. Until 
\cite{_BHPT:abelian_CY_} and \cite{_Kollar:Neron_},
N\'eron models were usually considered for abelian fibrations
over spectra of Dedekind rings,  but the complex
analytic context was developed only recently.

However, the relative Jacobian group scheme has been studied
in the hyperk\"ahler context since
\cite{_Markushevich:genus2_}
(see also \cite{_Arinkin_Fedorov_}
and \cite{_AbashevaRogov_}). Y.-J. Kim refined these
constructions to obtain a version of N\'eron
model for a Lagrangian fibration $\pi:\; X \to Y$.
 Unlike in the case of Bogomolov and
Koll\'ar, this group scheme is not realized as
an open subset of $X$. We construct a simplified version
of multi-dimensional N\'eron model, obtaining a group
or a torsor structure on an open subset of $X$, relative to $Y$.

The main statement of this paper appeared as a part of Theorem 1.6 
in \cite{_Kim:Neron_}. We give an independent proof of the
following statement (see also \cite[Theorem 1.6]{_Kim:Neron_}).

\hfill

%%%%%%%%%%%%%%%%%%%%%%%%%%%%%%%%%%%%%%%%%%%%%%%%%%%%%%%%%%%%
\theorem\label{_main_torsor_Intro_}
Let $X$ be a holomorphic symplectic manifold,
and $\pi:\; X \to Y$ a Lagrangian fibration, that
is, a proper holomorphic map which has holomorphic
Lagrangian fibers and admits a relatively ample line 
bundle.\footnote{Further in this paper, we call 
proper holomorphic maps admitting a relatively
ample line bundle {\bf projective}.} Assume that $Y$ is an open
subset of a projective manifold, and $X_i\subset X$
the set of all points where the differental
of $\pi$ has rank $\leq \dim Y-i$.\footnote{In 
\ref{_Diskr_i_codim_Claim_}, we show that $\codim \pi(X_i)\geq i$.}
 Denote by
$Y_0$ the set $Y\backslash \pi(X_2)$,
and let $N\subset X$ be $\pi^{-1}(Y_0) \backslash X_1$.
Then the natural projection $N \to Y_0$
is equipped with a unique structure of 
a torsor of a smooth group scheme over the base $Y_0$.

\proof
\ref{_Neron_any_dim_Theorem_}.
\endproof

\hfill

The proof of this result is based on the canonical model
construction of J. Koll\'ar (\cite{_Kollar:Neron_}). 
In \cite{_Kollar:Neron_},
the same statement is proven for any abelian fibration over a curve with
trivial relative canonical bundle (see \ref{_smooth_locus_group_Corollary_}). 
In assumptions of \ref{_smooth_locus_group_Corollary_},
the map $\pi:\; N \to Y_0$ is called {\bf the canonical model}
of the abelian fibration.

Our work generalizes 
\ref{_smooth_locus_group_Corollary_}
to the base of any dimension, assuming that the
abelian fibration is a Lagrangian fibration for
some holomorphic symplectic structure on its total space.
It is still unknown whether
a multi-dimensional version of Koll\'ar's result 
holds in a more general situation.

In \cite{_Kim:Neron_}, Y.-J. Kim has
significantly developed the construction 
of N\'eron models for holomorphic Lagrangian fibrations,
proving a version of \ref{_main_torsor_Intro_}
(and much more). Our approach is independent from his.

\hfill

Another novelty of this paper is our definition
of torsors, which significantly simplifies the exposition.
We define a torsor as a set $R$ with a map
$R\times R \times R \to R$, satisfying certain
algebraic conditions. The utility of this notion
is due to the fact that a compact complex torus
admits a unique structure of a holomorphic torsor.
Then, any abelian fibration can have at most
one structure of a torsor bundle over a group
bundle, and this torsor structure has at most
one extension to the special fibers. This approach
allows us to discuss the torsor analogue of the
N\'eron model for abelian fibrations which have
no sections.

The fibrations we consider are always proper morphisms, and have compact 
complex tori for general fibers.

%%%%%%%%%%%%%%%%%%%%%%%%%%%%%%%%%%%%%%%%%%%%%%%%%%%%%%%%%%%%
\section{N\'eron models}
%%%%%%%%%%%%%%%%%%%%%%%%%%%%%%%%%%%%%%%%%%%%%%%%%%%%%%%%%%%%

%%%%%%%%%%%%%%%%%%%%%%%%%%%%%%%%%%%%%%%%%%%%%%%%%%%%%%%%%%%%
\subsection{N\'eron models and canonical models}
\label{_Canonical_Nron_Subsection_}
%%%%%%%%%%%%%%%%%%%%%%%%%%%%%%%%%%%%%%%%%%%%%%%%%%%%%%%%%%%%

We introduce N\'eron models, starting from the
traditional, scheme-theoretic definition (see e. g. 
\cite{_BLR:Nron_models_,_BHPT:abelian_CY_}), and
proceed to their complex-geometric interpretation,
following J. Koll\'ar \cite{_Kollar:Neron_}.

\hfill

The notion of N\'eron model can be defined in bigger generality,
but for most practical purposes, people consider the N\'eron models
of abelian fibrations. In this context, ``abelian fibration''
is understood as a proper holomorphic map $\pi:\; X \arrow Y$ equipped
with a relatively ample bundle with
an abelian variety as a general fiber and a section. If $Y$ is quasiprojective,
$\pi$ defines an abelian variety over the field of rational
functions on $Y$, and this section represents its zero.

\hfill

Note that the abelian fibration, as defined above,
is equipped with a section, hence it has no multiple fibers.

\hfill

\definition
Let $R$ be a Dedekind domain, $K$ its fraction
field, and $Y$ a scheme over $K$.
For a scheme ${\cal X}$ over $R$
we denote the corresponding scheme
${\cal X}\times_{\Spec R}\Spec K$
over the fraction field by ${\cal X}_K$.
A {\bf N\'eron model} of $Y$ over $R$ is
a smooth, separated finite-type $R$-scheme 
${\cal Y}$ such that ${\cal Y}_K =Y$ 
and ${\cal Y}$ satisfies {\bf the N\'eron mapping
  property} (NMP): for any smooth $R$-scheme
${\cal X}$ and any $K$-morphism $f_K:\; {\cal X}_K\to {\cal Y}_K$,
$f_K$ can be uniquely extended to a morphism
$f:\; {\cal X}\arrow {\cal Y}$.

\hfill

\remark\label{_Nron_group_Remark_}
Uniqueness of a N\'eron model is clear.
As follows from \cite[Theorem 1.4/3]{_BLR:Nron_models_},
 the N\'eron model always exists for any abelian variety over $K$.
Also, the  N\'eron mapping property implies that the N\'eron model 
of an algebraic group over $K$ is a group scheme over $R$. 

\hfill

%This notion makes sense in the complex-geometric
%context. Let ${\cal Y}\stackrel \phi \arrow C$ 
%be an abelian fibration over a smooth complex curve $C$,
%and $K=k(C)$ its ring of rational functions.
%Since $C$ is covered by open sets which are
%isomorphic to spectra of Dedekind domains,
%it makes sense to define the N\'eron model of ${\cal Y}_K$
%as above. 
%
%\hfill

\remark The N\'eron model was defined over a base which is
a spectrum of a Dedekind domain. In the complex geometric
setting, this amounts to defining a N\'eron model
over the spectrum of a ring of germs of a 
smooth complex curve (which is a discrete valuation ring,
that is, a local Dedekind domain). This is a local version
of this notion, which can be globalized as follows. 
Define {\bf the N\'eron model} of an abelian fibration
$\pi:\; X \arrow C$
over a smooth complex curve $C$ as a complex variety
${\cal X}\stackrel \phi \arrow C$ bimeromorphic to $X$
which gives its N\'eron model if tensored by the ring
$\calo_{C,c}$  of germs of any point $c\in C$. 

\hfill

Note, however, that the Lagrangian fibrations we deal with
do not, generally speaking, admit sections, hence they
do not define abelian fibrations. If we restrict
the fibration to a sufficiently small open set in the
base, local sections exist, and we can always define
its N\'eron model.

\hfill

\theorem\label{_Nron_is_smooth_locus_Theorem_}
Let ${\cal X}\stackrel \phi \arrow C$  be an
abelian fibration over a smooth complex curve $C$.
Assume that ${\cal X}$ is smooth, algebraic, projective over $C$, and 
the relative canonical bundle of ${\cal X}/C$
is trivial. Then the N\'eron model of ${\cal X}_K$
is isomorphic to the smooth locus of $\phi$, which is
an open subset in ${\cal X}$.

\hfill

\proof In an earlier version of \cite{_BKV:Sections_}, we 
followed \cite{_BHPT:abelian_CY_}, Section 3 
and \cite[Proposition 5.1.8]{_HN:Motivic_zeta_},
referring to the notion of a ``Kulikov model''. However,
J\'anos Koll\'ar \cite{_Kollar:Neron_} pointed out that the definitions of the
Kulikov model in these papers are different, and,
for the proof, \cite{_BHPT:abelian_CY_}  refers to a claim in
\cite{_HN:Motivic_zeta_}, which does not exist. 

In the present version, 
we use a result of \cite{_Kollar:Neron_}.
Recall that {\bf a model} of ${\cal X}$
is a variety, fibered over $C$, bimeromorphic to
${\cal X}$ over $C$, and with generic fiber isomorphic
to that of ${\cal X}$. A {\bf minimal model}
of an abelian fibration is a model 
which has terminal singularities and numerically
trivial relative canonical bundle. Different flavors of this 
notion are called ``Kulikov model'' in
\cite{_BHPT:abelian_CY_} and \cite[Proposition 5.1.8]{_HN:Motivic_zeta_}.
By \cite[Corollary 2 (2.2a)]{_Kollar:Neron_}, 
the natural map from the smooth locus of 
a minimal model to the N\'eron model of ${\cal X}$ is an isomorphism.
\ref{_Nron_is_smooth_locus_Theorem_} is a special case of this theorem.
\endproof

\hfill

This result immediately brings the following corollary. It 
is essentially the only result from this subsection which we are going to use.

\hfill

%%%%%%%%%%%%%%%%%%%%%%%%%%%%%%%%%%%%%%%%%%%%%%%%%%%%%%%%%%%%
\corollary\label{_smooth_locus_group_Corollary_}
Let ${\cal X}\stackrel \phi \arrow C$  
be a projecive map to over a smooth complex curve $C$
with general fiber an abelian variety. Assume that $\phi$
admits holomorphic local sections everywhere in $C$.
Assume that ${\cal X}$ is smooth, algebraic, and 
the relative canonical bundle of ${\cal X}/C$
is pullback of a line bundle on $C$. Then the smooth locus of $\phi$
is equipped with a structure of an algebraic
group scheme over $C$.

\hfill

\proof
Locally in $C$, the bundle ${\cal X}/C$ admits holomorphic sections,
and therefore the smooth locus of $\phi$ is the N\'eron 
model of ${\cal X}$ (\ref{_Nron_is_smooth_locus_Theorem_}).
The N\'eron model is a group scheme by \ref{_Nron_group_Remark_}.
\endproof

%%%%%%%%%%%%%%%%%%%%%%%%%%%%%%%%%%%%%%%%%%%%%%%%%%%%%%%%%%%%
\subsection{Canonical model of Lagrangian fibrations}
%%%%%%%%%%%%%%%%%%%%%%%%%%%%%%%%%%%%%%%%%%%%%%%%%%%%%%%%%%%%

The canonical model can be applied to construct
a N\'eron model for abelian fibrations over a base
of dimension $>1$. This is an elementary alternative 
to a much more general strategy chosen by Y.-J. Kim
(\cite{_Kim:Neron_}).

\hfill

\definition
Let $f:\; X_0 \to Y$ be a smooth holomorphic submersion
(not necessarily proper), and $s:\; Y\to X_0$
a section. Consider a map $X_0\times_Y X_0\stackrel\mu \to X_0$
commuting with the projection to $Y$. We say that $(X_0,s, \mu)$
is {\bf a holomorphic group bundle} if the map $\mu$ is
submersive and defines a group structure on each fiber
$f^{-1}(y)$, for all $y\in Y$, and $s(y)$ is its identity element.

\hfill

%Recall that {\bf an abelian fibration} is 
% a projective, surjective map $f:\; X \to Y$,
%with general fiber an abelian variety. 
%Since $f$ is projective, its fibers are
%polarized, with the same K\"ahler class.
%Further on we will assume that the relative canonical
%bundle $K_f:= K_X\otimes f^* K_Y^{-1}$ is 
%isomorphic to a pullback of a line bundle on $Y$, 
%or, equivalently, that the canonical bundle $K_X$
%is isomorphic to a pullback of a line bundle on $Y$.
%This happens, in particular, when $X$ is 
%a Calabi-Yau manifold.
Further on, we will assume that 
$f:\; X \to Y$ is a Lagrangian fibration
on a hyperk\"ahler manifold, though
we expect that this construction might
work in bigger generality.

In this subsection, we construct
``the N\'eron model'' of such a morphism $f$,
 (called, following \cite{_Kollar:Neron_}, ``the canonical model'') 
when $f$ admits a holomorphic section
$s:\; Y \to X$. In this case, the canonical model is
a holomorphic group bundle.
In Subsection \ref{_torsor_Nron_Subsection_},
we generalize this notion to abelian fibrations
with no sections, in which case the canonical model 
is a smooth holomorphic abelian torsor bundle
(\ref{_holo_torsor_bundle_Definition_}).

\hfill

\definition
Let $f:\; X \to Y$ be an abelian fibration of
projective manifolds.  Let $X_0\subset X$ be the smooth
locus of $f$, that is,
the set of all points where $D\!f$ has maximal 
rank. Clearly, $X_0$ is a dense,
Zariski open subset of $X$. A smooth curve $C \subset Y$
is called {\bf N\'eron regular} if the smooth
locus of the map $f^{-1}(C) \stackrel f \to C$ is $X_0 \cap f^{-1}(C)$.

\hfill

We are going to show that a ``generic'' smooth curve
$C\subset Y$ is N\'eron regular.

\hfill

%%%%%%%%%%%%%%%%%%%%%%%%%%%%%%%%%%%%%%%%%%%%%%%%%%%%%%%%%%%%
\definition\label{_Discr_i_Definition_}
The {\bf discriminant} $\Discr\subset Y$ of $f$ is
the set of its critical values.
Let $X_i\subset X$ be the set of all points
$x\in X$ such that $D\!f\restrict x$ has rank $\dim Y-i$,
$\bar X_i$ its closure, 
and $\Discr_i := f(\bar X_{i+1})$. Clearly, $\Discr_0=\Discr$.

\hfill

%%%%%%%%%%%%%%%%%%%%%%%%%%%%%%%%%%%%%%%%%%%%%%%%
\claim\label{_Diskr_i_codim_Claim_}
For any $i$, we have $\codim_Y \Discr_i \geq i+1$.

\hfill

\proof
Since $\rk D\!f\restrict {\Discr_i} \leq \dim Y-i-1$,
Remmert rank theorem 
(\cite[page 36, Chapter 3, \S 3]{_MH_Schwartz_})
implies that $\dim f(\Discr_i) \leq \dim Y-i-1$.
\endproof

\hfill

%%%%%%%%%%%%%%%%%%%%%%%%%%%%%%%%%%%%%%%%%%%%%%%%
\claim
The complement $\Discr \backslash \Discr_1$ is smooth.

\hfill

\proof
By definition, 
$\Discr\backslash \Discr_1 = f(\bar X_1)\backslash f(\bar X_2)$,
where $\bar X_i$ denotes the closure. Also by definition,
$f$ has rank $\dim Y-1$ on $X_1$. Therefore, 
$\Discr\backslash \Discr_1$ is the set of regular values of
$f$ on $\bar X_1$.
\endproof

\hfill

%%%%%%%%%%%%%%%%%%%%%%%%%%%%%%%%%%%%%%%%%%%%%%%%%%%%%%%%%%%%
\proposition\label{_Nron_regular_Proposition_}
Let $f:\; X\to Y$ be an abelian fibration
and $C\subset Y\backslash \Discr_1$ a smooth curve 
such that the intersection $C\cap \Discr$
is transversal. Then $C$ is N\'eron regular.

\hfill

\pstep
Clearly, any $x\in f^{-1}(C)\cap X_0$ belongs to 
the smooth locus of $f^{-1}(C) \stackrel f \to C$.
Since $C\subset Y\backslash \Discr_1$, and
$\Discr_1= f(X_2)$, we have $f^{-1}(C)\subset X_1 \cup X_0$.
It remains to show that $f$ is not regular in 
any point $x\in  f^{-1}(C)\cap X_1$.

\hfill

{\bf Step 2:}
For any $x\in X_1$, the image $D\!f(T_x X)$
has codimension 1 in $TY$. By Sard lemma,
the map $X_1 \to \Discr$ is regular over
a general point $y\in \Discr$. Therefore, for any $x\in f^{-1}(y)\cap X_1$,
the image $D\!f(T_x X_1)$ had dimension $\dim Y-1$.
This implies that $D\!f(T_x X)= D\!f(T_x X_1)$, 
giving $D\!f(T_x X_1)\subset T_y Y$ 
for any $x\in f^{-1}(y)\cap X_1$.
Since the limit of tangent vectors
to a complex variety belongs to its
Zariski tangent space, and the Zariski tangent cone
of $\bar X_1$ is the closure of the Zariski tangent 
cone to $X_1$, this implies that
$D\!f(T_x X_1)\subset T_y \Discr$ for any $x\in \bar X_1$.

\hfill

{\bf Step 3:}
Let $x\in X_1$.
Since $D\!f(T_x X)$ has codimension 1 in $T_{f(x)} Y$
and $\dim \Discr \backslash \Discr_1 =\dim Y-1$,
the relation $D\!f(T_x X_1)\subset T_y \Discr$ 
implies that $D\!f(T_{\bar x} X)\subset T_y \Discr$
for any $\bar x\in X_1$.

\hfill

{\bf Step 4:}
Let $C\subset Y\backslash \Discr_1$ be a smooth curve 
such that the intersection $C\cap \Discr$
is transversal, and $x\in X_1 \cap f^{-1}(C)$.
Since $C$ is transversal to $\Discr$,
the tangent space $T_y C$ does not belong to $T_y \Discr$,
where $y=f(x)$. Since $D\!f(T_x X)\subset T_y \Discr$,
the differential $D\!f:\; T_x f^{-1}(C)\to C$
cannot be surjective, hence $x$ belongs to the
singular locus of the map $f:\; f^{-1}(C)\to C$.
\endproof

\hfill

\remark\label{_many_curves_Remark_}
Let ${\goth C}$ be a sufficiently big family of smooth curves
in $Y$, for instance, the set of all smooth intersections
of $Y\subset \C P^m$ with complex projective subspaces
of appropriate dimension. By Bertini's theorem,
a general curve $C\in {\goth C}$ is transversal to the smooth divisor 
$\Discr\backslash \Discr_1$ and does not intersect 
a codimension 2 subvariety $\Discr_1\subset Y$. 
Therefore, a general curve $C\in {\goth C}$
is N\'eron regular.

\hfill

Using \ref{_Nron_regular_Proposition_}, we 
observe that the smooth locus $f^{-1}(C)_0$ of
the map $f^{-1}(C)\stackrel f\to C$ 
coincides with the intersection
of the smooth locus $X_0$ of $f:\; X \to Y$
and $f^{-1}(C)$ for any N\'eron regular curve. 
\ref{_smooth_locus_group_Corollary_} implies that 
$f^{-1}(C)_0$ admits a natural group bundle structure if
$f$ is equipped with a section. It turns out that
this group bundle structure is actually obtained from
a group bundle structure on $X_0\cap f^{-1}(Y\backslash \Discr_1)$. 

\hfill

%%%%%%%%%%%%%%%%%%%%%%%%%%%%%%%%%%%%%%%%%%%%%%%%%%%%%%%%%%%%
\theorem\label{_canonical_model_group_Theorem_}
Let $X$ be a hyperk\"ahler manifold, 
$f:\; X \to Y$ a Lagrangian fibration,
and $X_0\subset X$ its smooth locus.
Denote by $J(f)$ the manifold $X_0\cap f^{-1}(Y\backslash \Discr_1)$.
Assume that $f$ admits 
a holomorphic section $s:\; Y \to X_0$.
Assume, moreover, that the relative canonical
bundle $K_f:= K_X\otimes f^* K_Y^{-1}$ is 
isomorphic to a pullback of a line bundle on $Y$.
Then there exists a unique group bundle structure
on $f:\; J(f)\to Y\backslash \Discr_1$ with the
following properties:
\begin{description}
\item[(i)] For each $y\in Y\backslash \Discr_1$,
$s(y)$ is the identity element of the fiber $f^{-1}\cap J(f)$.
\item[(ii)] For any N\'eron generic curve $C\subset Y$, 
the group structure on $J(f)$ is compatible with the
group structure defined
on $f^{-1}(C)\cap X_0$ by
\ref{_smooth_locus_group_Corollary_}.
\item[(iii)] 
The natural action of the associated group bundle $J(f)$ 
on itself extends holomorphically to $f^{-1}(Y\backslash \Discr_1)$.
\end{description}

\pstep
The uniqueness of such a group structure is clear.
Indeed, holomorphic group structures on a compact complex torus
are in bijective correspondence with the choice of an identity
element. Since a general fiber of $f$ is a compact torus,
the group structure on the special fibers
of $f:\; J(f)\to Y\backslash \Discr_1$ is 
unique by continuity.

\hfill

{\bf Step 2:}
To prove the existence of the group 
bundle structure on $f:\; J(f)\to Y\backslash \Discr_1$,
consider a special fiber $f^{-1}(y)\subset J(f)$.
Then $y\in \Discr\backslash \Discr_1$. Since
$Y$ is projective, it admits a family ${\goth C}$ of
N\'eron regular curves  satisfying the
assumptions of \ref{_many_curves_Remark_}.
For a general curve $C\in {\goth C}$ containing $y$, the intersection
$J(f)\cap f^{-1}(C)$ admits a group bundle
structure. This group bundle structure,
in turn, defines a group structure on $f^{-1}(y) \cap J(f)$.
To finish the proof of 
\ref{_canonical_model_group_Theorem_},
it remains to show that this group structure
is independent of the choice of $C \in {\goth C}$.

\hfill

{\bf Step 3:}
Let $u_1, ..., u_n\in \calo_Y$ be a collection
of holomorphic functions defined in 
a neighbourhood $U\subset Y$ of $y\in Y$, such that
their differentials generate $\Omega^1U$,
and let $\xi_1, ..., \xi_n\in T(f^{-1}U)$ the
corresponding holomorphic Hamiltionian vector fields.
They are defined by the relation $i_{\xi_j}\Omega=du_j$,
where $i_{\xi_j}$ denotes the contraction and $\Omega$
the holomorphic symplectic form. Since $du_i$ vanish
on the tangent space to the fibers of $f$, which are Lagrangian,
the vector fields $\xi_i$ are tangent to the fibers of $f$.
They are also non-degenerate and linearly independent
on the smooth locus of $f$. Therefore, they trivialize
the tangent bundle $Tf^{-1}(Y)\cap X_0$
to the smooth part of the special fiber $f^{-1}(y)$.

The Lie group structure on the fibers of $f$
can be obtained by taking the exponent
of its Lie algebra which is generated by $\xi_i$.
This exponent is well defined in the special
fiber as well. By continuity of the group structure
defined by the N\'eron model of $f^{-1}(C)$, it is also obtained
by taking exponents of $\langle \xi_1, ..., \xi_n\rangle$.
Therefore, it is independent of the choice of $C$.

\hfill

{\bf Step 4:}
To prove (iii), we notice that the Hamiltonian vector
fields defined in Step 3 extend to holomorphic
vector fields on the complex manifold $f^{-1}U$;
the corresponding group action coincides with the
action of $J(f)$ on $X_0$.
\endproof

\hfill

%%%%%%%%%%%%%%%%%%%%%%%%%%%%%%%%%%%%%%%%%%%%%%%%%%%%%%%%%%%%
\remark\label{_open_subset_Neron_Remark_}
The statement of \ref{_canonical_model_group_Theorem_}
can be formulated in a less restrictive situation,
without assuming that $X$ is hyperk\"ahler. This theorem
is true for any projective Lagrangian fibration over a base
which admits sufficiently big families of smooth curves.
For instance, it is true for $X = \pi^{-1}(U)$,
where $\pi:\; M \to B$ is a Lagrangian fibration on
a hyperk\"ahler manifold and $U\subset B$ an open subset.

%%%%%%%%%%%%%%%%%%%%%%%%%%%%%%%%%%%%%%%%%%%%%%%%%%%%%%%%%%%%
\subsection{Torsor bundles}
\label{_torsor_Nron_Subsection_}
%%%%%%%%%%%%%%%%%%%%%%%%%%%%%%%%%%%%%%%%%%%%%%%%%%%%%%%%%%%%

By definition, a torsor over an abelian group
$G$ is a set where $G$ acts freely and transitively.
However, the torsor structure can be defined
without referring to an {\em a priori} given group structure as follows.

\hfill

\definition
{\bf An abelian torsor} is a set $S$ equipped
with a map $\mu:\; S\times S \times S \to S$ such that
\begin{description}
\item[(i)] for
any $e\in S$, the map $x \cdot y := \mu(x, y, e)$
defines a commutative and associative product on $S$
with $e$ an identity element 
\item[(ii)] 
 for any fixed $y \in S$,
the map $x\mapsto x\cdot y$ is bijective.
\end{description}

\remark
From the axioms of the abelian
torsor, it is clear that $\cdot$ defines a structure 
of an abelian group on $S$.

\hfill

\remark
The identity element axiom means that 
$\mu(x,y,y)=x$ for any $x$ and $y$. Given a group
structure on $S$, the torsor operation
$\mu$ can be defined  as $\mu(x, y,z)= x+y-z$.

\hfill

This language is best suited to 
define N\'eron models for Lagrangian fibrations
which have no sections.

\hfill

%%%%%%%%%%%%%%%%%%%%%%%%%%%%%%%%%%%%%%%%%%%%%%%%%%%%%%%%%%%%
\definition\label{_holo_torsor_bundle_Definition_}
Let $f:\; X \to Y$ be a holomorphic submersion,
not necessarily proper, and 
$\mu:\; X\times_Y X\times_Y X\to X$
a holomorphic map commuting with the projection to $Y$,
which defines an abelian torsor structure in each
fiber of $f$. Then $(f, \mu)$ is called {\bf a smooth
holomorphic abelian torsor bundle}
over $Y$, or simply {\bf a torsor bundle}.

\hfill

\remark
Sincee $f$ is submersive, $f:\; X \to Y$ 
admits local holomorphic sections by the implicit map
theorem. Given a local section $s:\; V \to U$,
we obtain a group bundle structure on $f:\; f^{-1}(V)\to V$.
Clearly, this group bundle is independent of the
choice of $s$, hence we might glue the local group
bundles obtained this way to a group bundle $J(f)$ over
$Y$. This group bundle acts on $X$ holomorphically, and
this action is free and transitive on all fibers.
In other words, $X$  is a torsor bundle
over the group bundle $J(f)$.

\hfill

We can reformulate \ref{_Nron_is_smooth_locus_Theorem_}
using this language.

\hfill

%%%%%%%%%%%%%%%%%%%%%%%%%%%%%%%%%%%%%%%%%%%%%%%%%%%%%%%%%%%%
\proposition\label{_torsor_bundle_from_Nron_Proposition_}
Let ${\cal X}\stackrel \phi \arrow C$  be an
abelian fibration over a smooth complex curve $C$.
Assume that ${\cal X}$ is smooth, algebraic, projective over $C$, and 
the relative canonical bundle of ${\cal X}/C$
is trivial. Then the smooth locus of $\phi$
admits a natural structure of an abelian
torsor bundle, which is uniquely determined
by $\phi$.

\hfill

\proof
Existence of the torsor structure
is already explained in the proof of \ref{_Nron_is_smooth_locus_Theorem_};
it follows from \cite[Corollary 2 (2.2a)]{_Kollar:Neron_}.
The uniqueness of the torsor structure on 
an abelian variety is clear because it is 
the unique torsor structure compatible with 
the flat coordinates, which are intrinsically 
defined by the holomorphic differentials.
Since the torsor structure map
$\mu:\; {\cal X}\times_C {\cal X}\times_C {\cal X}\to {\cal X}$
is continuous, it is uniquely extended
to the special fibers.
\endproof

\hfill

This observation brings the following remarkable theorem.

\hfill

%%%%%%%%%%%%%%%%%%%%%%%%%%%%%%%%%%%%%%%%%%%%%%%%%%%%%%%%%%%%
\theorem\label{_Neron_any_dim_Theorem_}
Let $f:\; X\to Y$ be
a Lagrangian fibration on a hyperk\"ahler manifold,
$\Discr$ its discriminant, and $\Discr_1 \subset \Discr$
the locus of $\Discr$ defined in \ref{_Discr_i_Definition_}.
Denote by $Y_0\subset Y$ the set of all
points $y\in Y\backslash \Discr_1$ such that $f$ has local
holomorphic sections over a neighbourhood of $y$ with values
in the smooth locus of $f$. Denote by $X_0'$ 
the intersection of the smooth locus of $f$ with
$f^{-1}(Y_0)$, and let 
$Y'_0:= Y\backslash \Discr_1$. Then $f:\; X_0'\to Y'_0$
admits a unique structure of an abelian torsor bundle,
compatible with the torsor bundle structures
on $f^{-1}(C)$ obtained in 
\ref{_torsor_bundle_from_Nron_Proposition_},
for any N\'eron regular curves $C \subset Y_0$.
Moreover, the natural action of the associated group bundle $J(f)$ 
extends holomorphically to $f^{-1}(Y_0)$.

\hfill

\proof
By definition of $X_0$, the map $f:\; X_0'\to Y'_0$
 admits local sections. Let $U\subset Y'_0$
be an open subset and $s:\; U \to X_0'$ a holomorphic section.
Then  \ref{_canonical_model_group_Theorem_} and 
\ref{_open_subset_Neron_Remark_}
produce a unique structure of a group bundle on 
the fibration $f^{-1}(U)\cap X'_0\to U$,
such that $s(y)$ gives the identity element
in each fiber $f^{-1}(y)$. A different choice $s_1:\; \; U \to X_0'$ 
of a section produces a different group structure,
but the corresponding group bundles are isomorphic.
The isomorphism map is given by a parallel translation
on each fiber, $x\mapsto x- s(f(x))+ s_1(f(x))$.
Therefore, the torsor structure on the fibration
$f^{-1}(U) \to U$ is independent from the choice of a section $s$.
Gluing together the structure maps
$f^{-1}(U)\times_U f^{-1}(U)\times_Uf^{-1}(U)\to U$,
we obtain a structure of a torsor 
$X'_0\times_{Y'_0}  X'_0\times_{Y'_0}X'_0\to Y'_0$
extending the standard torsor structure
on the general fibers of $f:\; X'_0 \to Y'_0$.
\endproof

\hfill

{\bf Acknowledgments.} 
We are grateful to F. Bogomolov and J. Koll\'ar
for their interest and many enlightening discussions,
and to Y.-J. Kim for his advice and his interest.

\hfill

\hfill

\small
\noindent {\sc Ljudmila Kamenova\\
Department of Mathematics, 3-115 \\
Stony Brook University \\
Stony Brook, NY 11794-3651, USA,} \\
\tt kamenova@math.stonybrook.edu
\\

\noindent {\sc Misha Verbitsky\\
            {\sc Instituto Nacional de Matem\'atica Pura e
              Aplicada (IMPA) \\ Estrada Dona Castorina, 110\\
Jardim Bot\^anico, CEP 22460-320\\
Rio de Janeiro, RJ - Brasil}\\
\tt verbit@impa.br
}

\end{document}